\hsize=16 true cm
\vsize=20.5 true cm
  \baselineskip=15pt
\vglue 1.3 cm
\count1=\number\year
\advance\count1 by -1900
\def\abc{\number\month/\number\day/\number\count1}
\footline={\hss\tenrm -\folio-\ \ \abc \hss}
        


  \def\qed{$\rlap{$\sqcap$}\sqcup$}


    \font\tenmsb=msbm10              \font\sevenmsb=msbm7
\newfam\msbfam
      \textfont\msbfam=\tenmsb
      \scriptfont\msbfam=\sevenmsb
\def\Bbb#1{{\fam\msbfam #1}}



\def\move-in{\parshape=1.75true in 5true in}


\hyphenation {Castel-nuovo}
           

\def\PP#1{{\Bbb P}^{#1}}
\def\ref#1{[{\bf #1}] }  
\def\prim{^\prime } 
\def\fnote#1{\footnote{$^{#1}$}}

\def\Prop#1{\noindent {\bf Proposition #1:} }


\def\Coroll#1{\noindent {\bf Corollary #1:} }


\magnification = 1000
\vglue 10 pt
\centerline{{\bf BLOWING UP DETERMINANTAL SPACE CURVES }}
\centerline{{\bf AND VARIETIES OF CODIM 2}~\fnote{(*)}{This work has
been supported by MURST
funds and by the CNR group "Rami analitici e sistemi lineari".}}
\bigskip
\centerline {A.Gimigliano, A.Lorenzini}
\vskip 1cm
\noindent{\bf Introduction.}
\medskip
Let $k$ be an algebrically closed
field with char$k = 0$ and let $C\subseteq \PP 3 = \PP 3 _{k}$ be
a curve (i.e. a 1-dimensional, smooth, irreducible scheme).
\par
We want to study the scheme $X_C$ which is the
blow-up of $\PP 3$ along $C$. Let $E$ be the exceptional divisor of
the blow-up, and $H$ the strict transform of a generic hyperplane, then
${\rm Pic}X_C = {\bf Z}<H,E>$, and the kind of questions we would like
to address are:
\par
\medskip
\noindent 1)  For which $p$ is the linear system $\vert pH-E\vert$
very ample on $X_C$?
\par
\medskip
\noindent 2)  When is $\vert pH-E\vert$ normally generated ?
 \par
\medskip
\noindent 3) Let  ${\cal O}_{X_C}(pH-E)$ be generated by its global sections
(g.b.g.s. for short), then it defines a morphism
$\phi_p:X_C \rightarrow \PP N$, where
$N+1 = \dim H^0({\cal O}_{X_C}(pH-E))$.
What can we say about the ideal of $Y_{C,p}=\phi_p(X_C)$?
More specifically, we try to get information about the generators and
the resolution of its ideal  from the data known about $C$.
\par
\medskip
Note that by $\vert pH-E\vert$ {\it normally generated} we mean
that ${\cal O}_{X_C}(pH-E)$ is very ample and $Y_{C,p}$ is
arithmetically Cohen-Macaulay (a.C.M. for short), i.e. its coordinate
ring is Cohen-Macaulay.
\par
\medskip
We will restrict ourselves to the case when $C$ is projectively normal.
Let  $I_C \subseteq R = k[w_0,w_1,w_2,w_3]$
be its  homogeneous ideal; by the
Hilbert-Burch Theorem
we know that $I_C$ is generated by the maximal minors of
a $\rho\times (\rho+1)$ matrix :
$$\qquad M = \pmatrix{M_{11}&M_{12}&\quad ...\quad &M_{1\rho +1}\cr
M_{21}&M_{22}&\quad ...\quad &M_{2\rho +1}
\cr \vdots & \vdots&&\vdots \cr
M_{\rho 1}&M_{\rho 2}&\quad ...\quad &M_{\rho \rho +1}}$$
where the $M_{ij}$'s are forms in $R$ and
$\deg M_{ij} = \deg M_{i1} + \deg M_{1j} - \deg M_{11}$.
The degree matrix, $\partial M$, is determined by the degrees of
a minimal set of generators of $I_C$ or, equivalently in this case,
 by the graded Betti numbers of a minimal resolution of $I_C$.
\par
A resolution of the ideal sheaf ${\cal I}_C$ is of the type:
$$0 \rightarrow \oplus _{j=1}^u{\cal O}_{\PP 3}(-n_j) \rightarrow
\oplus _{i=1}^v{\cal O}_{\PP 3}(-d_i) \rightarrow {\cal I}_C \rightarrow 0.$$
\par
\medskip
A first answer to question 1) above is given by  [{\bf Co}, Theorem 1]:
let $V\subseteq \PP n$ be a smooth irreducible projective scheme
and $\lambda$  such that
$V$ is scheme-theoretically generated in degrees $\leq \lambda$, then
 (with obvious notation) $\vert pH-E\vert$ is very ample on the blow up $X_V$
of $\PP n$ along $V$ for all $p\geq \lambda +1$.
\par
\medskip
In particular, let $V$ be a.C.M. and consider the  invariant
$$
\sigma = \sigma(V) = min \{t\vert \Delta^{n-1}H(V,t)=0\},
$$
\noindent where $H(V,t)$ is the Hilbert function of $V$,
$\Delta^{i}H(V,t)=\Delta^{i-1}H(V,t)-\Delta^{i-1}H(V,t-1)$ and
 $\Delta^{0}H(V,t)=H(V,t)$.
 \par
\medskip
It is well known that the homogeneous ideal of $V$, $I_V$,
is generated in degrees $\leq \sigma$
 (see e.g. \ref {L.1, Theorem 2}),
hence $\vert pH-E\vert$ is always very ample for $p\geq \sigma +1$.
\par
\medskip
Several examples and known facts lead to the following conjecture:
\par
\bigskip
\noindent {\bf Conjecture}: {\it Let} $C\subseteq \PP 3$ {\it
be a projectively normal curve.  Then}
$\vert \sigma H-E\vert$ {\it is very ample on } $X_C$ {\it if and only if
C has no} $\sigma${\it -secant lines}. {\it  Moreover, if this is the case,}
$\vert \sigma H-E\vert$ {\it  is normally generated.}
\par
\bigskip
Note that
${\cal O}_{X_C}( \sigma H-E)$ is generated by its global section since
$I_C$ is generated in degrees $\leq \sigma$.
\par
The bound $p\geq \sigma$ is sharp with no other specification on
$C$, since $\vert \sigma H-E\vert$ is not very ample when $C$ has a
$\sigma$-secant $L$ (the divisors in $\vert \sigma H-E\vert$
will not intersect the strict transform of $L$).
Examples of curves $C$
possessing $\sigma$-secants are well known:
\par
- rational normal cubics curves ($\sigma = 2$), the most trivial example;
\par
- projectively normal sextic curves of genus 3
(here $\sigma = 3$ and they have infinitely many 3-secants, see e.g. \ref R);
\par
- quintic curves of genus two of type (3,2) on a quadric surface
( $\sigma = 3$).
\medskip
The whole problem can be generalized in a very natural way to the blow-up,
$X_V$, of $\PP n$ along a codimension 2 (smooth, irreducible)
variety $V\subseteq \PP n$, $n\geq 3$. If we require that
$V$ is a.C.M. we still have that the ideal $I_V$ of $V$ is
generated
in degrees $\leq \sigma$ and, by Hilbert-Burch theorem, determinantal.
\par
Notice that the analog of the Conjecture above for points in $\PP 2$ is true
(see \ref {D-G}), but
we will not consider this case here, since it has been
studied by many authors (e.g. see the surveys \ref {Gi.2}, \ref {Gi.3}).
\par
Note also that, when $\rho \geq 2$, the locus $\Sigma $
where the rank of the matrix $M$ drops twice is singular for $V$
and its codimension is
$\leq   6$, hence,
for $n\geq 6$, $\Sigma$ is non-empty and $V$ has to be singular.
\par
Thus, for $n\geq 6$, only $\rho =1$ is allowed in order to have $V$ smooth,
i.e. $V$ has to be a complete intersection,
so its ideal is quite simple and if it is generated by two forms of
degrees $d_1,d_2$, say with $d_1\geq d_2$,
we have that $\vert pH-E\vert$ is very ample if
$p\geq d_2+1$ (e.g. by \ref{Co}). We will not treat that case here, see e.g.
\ref {S-T-V} and \ref {G-G-H} for the case  $n=2$.
\par
\medskip
The plan of the paper is the following:
\par
\medskip
\noindent \S 1: In the general case, we show that $\vert \sigma H-E\vert$ is
 very ample on
$X_V-E$ (i.e. $(I_V)_\sigma$ gives a very ample linear system on $\PP n - V$).
\par
\medskip
\noindent \S 2: We prove our Conjecture  when
the entries of the matrix $M$ are linear forms. In this case a quite
complete description
of the ideal of  $Y_{C,\sigma}$
(which turns out to be determinantal itself)
and of its resolution can be given.
\par
Furthermore, we extend our result to the case
of codim 2 subschemes of $\PP n$, $n=4,5$, which do not
contain lines and
we give numerical conditions that describe
when the generic determinantal scheme of the type under
consideration does
contain lines or
possess $\sigma$-secants.
\par
\medskip
\noindent \S 3: For a given degree $d \geq 3$, we will consider
projectively normal curves $C$ of degree $d$
 with a "generic" resolution  (these curves have minimal genus for that
degree).
\par
In this case we are able to find the ideal generation of
$Y_{C,\sigma}$.
\par
\medskip
\noindent \S 4:  We give some examples of the construction seen above in the
 case of Fano varieties of dimensions 3,4,5.
\par
\bigskip
We would like to thank Marvi Catalisano and Tony Geramita for some
useful talks.
\par
\vskip 1cm
\noindent {\bf 1. Very ampleness on $\PP n -V$}.
\medskip
Let $V\subseteq \PP n$ be a smooth, irreducible, a.C.M. scheme of codimension
2 and let $\sigma$, $X_V$ be as in the introduction. The
first step in studying the very ampleness of
$\vert \sigma H-E\vert$ is to show that
$\phi _\sigma$ is generically 1:1; more precisely we will show:
\bigskip
\noindent {\bf Proposition 1.1:} {\it If there are no } $\sigma$-{\it secant
lines to V, the linear system}
$\vert \sigma H-E\vert$ {\it separates points and tangent vectors
on} $X_V-E$.
\par
\bigskip
Note that this is equivalent to saying that the sections of
$H^0(\PP n,{\cal I}_{V}(\sigma)) \cong (I_{V})_\sigma$
separate points and tangent vectors on $\PP n - V$.
\par
\medskip
\noindent {\it Proof:} What we have to show is that given any
0-dimensional scheme $T\subseteq \PP n - V$ of degree 2, we have that
the restriction map:
\par
\medskip
\noindent$(\dagger ) \qquad \qquad \qquad \qquad \qquad
H^0(\PP n,{\cal I}_{V}(\sigma)) \rightarrow
H^0(\PP n,{\cal O}_{T}(\sigma))$
\par
\medskip
\noindent is surjective.
\par
Let $L$ be the line defined by T; we will be done if
we show that the restriction of $I_\sigma = (I_{V})_\sigma$ to $L$
defines a very ample linear system on $L - (V\cap L)$ if and only if
$L$ is not a $\sigma$-secant for $V$.
\par
Let $I_L$ be the homogeneous ideal of $L$ and consider the exact sequence:
$$\matrix{&&&&&{}_\alpha&\cr 0& \rightarrow& (I_V\cap I_L)_\sigma &
\rightarrow & (I_V)_\sigma &
\rightarrow & (I_{V\cap L,L})_\sigma}$$
\par
\medskip
Since $V\cap L$ is given by $k$ points (not necessarely distinct)
on $L$, and $0\leq k \leq \sigma -1$, we have
$\dim (I_{V\cap L,L})_\sigma = \sigma -k+1\geq 2$,
hence it defines, out of $V\cap L$,
a $g_{\sigma -k}^{\sigma^2-k}$ on $L\cong\PP 1$, which is very ample
on $L - (V\cap L)$ because it does not have base points (out of $V\cap L$)
since $I_{V\cap L,L}$ is generated in degree $k<\sigma$.
\par
So we will be done if $\alpha $ is surjective.
\par
\medskip
\noindent {\bf Fact:} Let $\Pi$ be a generic plane containing
$L$ and let $Z=V\cap \Pi$,
then $Z$ is a 0-dimensional scheme
of degree $d = \deg V$ in  $\Pi \cong \PP 2$ and the homogeneous ideal of
$Z$ in $\Pi$ is determinantal (with the same degree matrix as $V$).
\par
\medskip
In order to check that such $\Pi$ exists one
has to show that there are
$n-2$ hyperplanes $H_1,...,H_{n-2}$ such that
$V_i = V\cap H_1\cap ...\cap H_{i}$ has dimension
$n-2-i$, $i=1,..,n-2$ (then we will choose $\Pi = H_1\cap...\cap H_{n-2}$).
\par
Consider the linear system ${\bf W} \subseteq H^0(V,{\cal O}_V(1))$ which
is cut on $V$ by the hyperplanes through $L$; since  $\dim (L\cap V) = 0$,
${\bf W}$ has no fixed components of dimension $\geq 1$.  Thus, if we
start with any $H_1$ through $L$, and
$D_1 = V\cap H_1 \in {\bf W}$,
 not all $D_2 \in {\bf W}$ will contain a component of $D_1$, hence
$D_2\in {\bf W}$ can be chosen such that $\dim (D_1\cap D_2) = n-4$ ($D_1$ and
$D_1\cap D_2$ are of pure dimension, since $V$ is a.C.M. and we cut
each time with a non-zero divisor).
\par
Iterating this procedure, at each cut the dimension drops by 1
untill  $D_1\cap ... \cap D_{n-2}$ has dimension zero, as required.
\par
\medskip
We have another exact sequence:
$$\matrix{&&&&&{}_{\alpha\prim}&\cr 0& \rightarrow& (I_Z\cap I_L)_\sigma &
\rightarrow & (I_Z)_\sigma &
\rightarrow & (I_{Z\cap L,L})_\sigma}$$
\medskip
\noindent and we will be done if $\alpha\prim$ is surjective, i.e. if
\par
\medskip
$\dim (I_Z\cap I_L)_\sigma =
\dim (I_{Z})_\sigma - \dim (I_{Z\cap L,L})_\sigma =
{\sigma + 2\choose 2} - d - (\sigma + 1 - k) = {\sigma + 1\choose 2} - d + k$.
\par
\medskip
Let  $Z\prim = Z - (Z\cap L)$ as schemes (i.e. $Z\prim$ is the
subscheme of
$\Pi$ defined by $I_Z:I_L$); note that
$\dim (I_Z\cap I_L)_\sigma = \dim (I_{Z\prim})_{\sigma -1}$,
since all the forms in $(I_Z\cap I_L)_\sigma$ are the product of a linear
form (defining $L$ on $\Pi$) and a form of degree $\sigma - 1$ containing
$Z\prim$. Thus all we need to show is that:
$$\dim (I_{Z\prim})_{\sigma -1} =  {\sigma + 1\choose 2} - d + k $$
\medskip
\noindent
i.e. that $Z\prim $ imposes independent conditions to the curves
of degree $\sigma -1$ (in fact $\deg Z\prim = d-k$), but this is
trivial, since $Z\prim\subseteq Z$, and $Z$ already does (it
follows from the
definition of $\sigma $).
\par
{\rightline {\qed}}
\bigskip
\noindent {\bf Remark 1.2:}\par
Note that the above proof works on every line $L\subseteq \PP n$, $L
 \not\subseteq V$,
hence
it also shows that $\vert \sigma H - E\vert $ separates points and
tangent vectors lying on the strict transform of such a line.
\par
\bigskip
Now let $\pi : X_V \rightarrow \PP n$ be the blowing up,
$P\in V$, and $\pi^{-1}(P)=F\subseteq E$ an $(n-2)$-dimensional
linear space in the ruling of $E$.
\par
\bigskip
\noindent{\bf Proposition 1.3:} {\it The linear system}
$\vert \sigma H-E\vert$ {\it is very ample on F}.
\par
\medskip
\noindent {\it Proof:}
This is quite an elementary fact, see e.g.
\ref{Co, \S 1}: we have to show that points $P\in F$ and tangent
vectors $t \in T_P(F)$ can be separated by divisors
$D \in \vert \sigma H-E\vert$.
The passage of $D$ through $P$
and its tangency to $t$ are equivalent to the fact that
$\pi(D)$ has certain tangent vectors at $\pi(P)$, and, by \ref {Co,
Lemma 1.5}, for every hyperplane $\Lambda \in T_{\pi (P)}(\PP n)$, one
can find $D \in \vert \sigma H-E\vert$ such that
$T_{\pi (P)}(\pi(D)) = \Lambda$, thus one can separate
tangent vectors at $\pi(P)$.
\par
\rightline {\qed}
\bigskip
\noindent {\bf Remark 1.4:}\par
Note that the above proof works also to show that $\forall P \in E$,
vectors $t\in T_P(X_V)$ can be separated by divisors
$D \in \vert \sigma H-E\vert$.
\par
\vskip 1cm
\noindent {\bf 2. Varieties defined by matrices of linear forms.}
\medskip
Let $V$ be as in \S 1 and let its defining matrix $M$ be such that
all its entries are linear forms.  In this case we have $\rho = \sigma$ and
it is quite easy to check that:
\bigskip
\Prop {2.1} {\it The variety} $Y_{V,\sigma}=\phi_\sigma (X_V)$ {\it is
a.C.M.}
\medskip
\noindent {\it Sketch of proof :}
 Let $M = (L_{ij})$,
 where $L_{ij}=\sum_{k=0} ^n \delta_{ij}^kw_k $. The homogeneous ideal of
$Y_{V,\sigma}$, $I_Y \subseteq k[x_1,...,x_{\sigma +1}]$,
is generated by the maximal minors
of a $((n+1)\times \sigma )$-matrix $N$ of linear forms in the $x_j$'s,
and since  ht$I =$ cod $Y_{V,\sigma} = \sigma - n$, $I_Y$ is perfect
and the matrix $N$ has entries
$N_{ik}=\sum_{j=1} ^{\sigma +1} \delta_{ij}^kx_j$
(e.g. see \ref {E-N}, \ref R , \ref {Gi.1}).
\par
\rightline {\qed}
\bigskip
Now it is quite immediate to give a resolution of the ideal $I_Y$;
\par
\bigskip
\Coroll {2.2} {\it A minimal resolution of} $I_Y$ {\it is given by :}
$$0 \rightarrow \oplus ^{{\sigma -1\choose n}}S(-\sigma) \rightarrow
... \rightarrow  \oplus^{{2+i\choose n}{\sigma \choose n+i}}S(-n-i)
\rightarrow
... \rightarrow \oplus ^{{\sigma\choose n+1}}S(-(n+1)) \rightarrow
I_Y \rightarrow 0.$$
{\it where}$\quad i=1,...,\sigma-n$ {\it and} $S = k[x_0,...,x_\sigma]$
{\it is the coordinate ring of} $\PP \sigma$.
\par
\medskip
\noindent {\it Proof:}
>From the proof of 2.1, we have that $I_Y$ is generated by the
maximal minors of the matrix $N$, so
the resolution above is
given by an Eagon-Northcott complex (see \ref{E-N}).
\par
\line{\hfill  \qed }
\bigskip
Now let us go back to the case of curves in $\PP 3$:
in particular we will consider the case when the curve $C$ is defined
by a matrix $M$ of linear forms; in this case
we are able to prove the Conjecture, namely we have the following:
\par
\bigskip
\noindent {\bf Theorem 2.3:}
{\it Let C be a projectively normal space curve
and let
$I_C$ be generated by the maximal minors of a matrix of linear forms.
Then} $\vert \sigma H-E\vert$ {\it is very ample on} $X_C$ {\it if
and only if C does not have any} $\sigma${\it -secant.}
\par
\medskip
\noindent{\it Proof:} We only have to show that the absence of
$\sigma$-secants
implies the very ampleness of $\vert \sigma H-E\vert$, since we already
noticed
in the introduction that this is a necessary condition for very ampleness.
\par
In the case when $C$ is defined by
a matrix $M$ of linear forms in $k[w_0,...,w_3]$ we know that $M$ is a
$\sigma\times(\sigma +1)$-matrix, that
$d = \deg C = {\sigma +1 \choose 2}$ and
the genus of $C$ is
$ g = 2{\sigma +1 \choose 3} - {\sigma +1 \choose 2} + 1$
(e.g. see \ref R).
\par
We will indicate:
$$\qquad M = \pmatrix{L_{11}&L_{12}&\quad ...\quad &L_{1\sigma +1}\cr
L_{21}&L_{22}&\quad ...\quad &L_{2\sigma +1}
\cr \vdots & \vdots&&\vdots \cr
L_{\sigma 1}&L_{\sigma 2}&\quad ...\quad &L_{\sigma \sigma +1}}\qquad
{\rm where}\qquad
L_{ij}=\sum_{k=0} ^3 \delta_{ij}^kw_k .$$
Let $F_{\beta} = (-1)^{\beta} F\prim_{\beta}$, where $F\prim _\beta$ is the
minor of $M$ obtained by erasing its $\beta^{th}$ column and
consider the map $\phi :\PP 3 - C \rightarrow \PP \sigma $ given by
$\phi (P) = [F_1(P):...:F_{\sigma +1}(P)]$ (we have that
$\phi_\sigma = \pi \circ \phi$ on $X_C - E$).
\par
Note that we surely have $\sigma \geq 3$, otherwise $C$ would
possess 1- or 2-secants.
\par
\medskip
We have to show that $\vert \sigma H-E\vert$ is very ample, i.e.
that given any 0-dimensional scheme $T \subseteq X_C$, with
$\deg T = 2$, the restriction map:
$$H^0(X_C,{\cal O}_{X_C}(\sigma H-E)) \rightarrow
H^0(X_C,{\cal O}_{T}(\sigma H-E)) = H^0(X_C,{\cal O}_T))$$
\noindent is surjective.
\par
We are already done if either:
\par
\medskip
\noindent - \quad $T\subseteq F \subseteq E$, where
$F$ is a fiber of $\pi$,  by Proposition 1.3;
\par
\noindent - \quad $T\subseteq X_C-E$,  by Proposition 1.1;
\par
\noindent - \quad Supp $T = P \in E$,  by Remark 1.4.
\par
\medskip
Hence we are left with the following possibilities:
\par
\medskip
\noindent {\bf Case 1)}:\quad $T=Q_1\cup Q_2 \subseteq E$, $Q_1\neq Q_2$,
and, if $P_i=\pi (Q_i)$, $P_1\neq P_2$;
\par
\noindent {\bf Case 2)}:\quad $T=Q_1 \cup Q_2$,  $Q_1\in E$,  $Q_2\notin E$.
\par
\medskip
\noindent {\bf Case 1)}:
There exist vectors
$t_i \in T_{P_i}(\PP 3)-T_{P_i}(C)$, $i=1,2$  such that
$\forall D \in \vert\sigma H-E\vert$, $Q_i \in D$
if and only if $t_i \in T_{P_i}(\pi (D))$; see e.g. \ref {Co}, \S 1.
\par
\medskip
Let $L_i\subseteq \PP 3$ be the line determined by $t_i$; if $L_1=L_2=L$,
then $T$ is contained in the strict transform of $L$, and we are done by
Remark 1.2 .
\par
\medskip
Let $L_1\neq L_2$, and suppose that $L_1\cap L_2=\emptyset$,
 then, with a linear change of coordinates,
we can assume that $P_1 =(1:0:0:0)$,
$P_2 =(0:0:0:1)$, $L_1 = \{w_2=w_3=0\}$, $L_2=\{w_0=w_1=0\}$.
\par
\medskip
What we have to check is that for surfaces $S=\{f=0\}$, with
$f \in (I_C)_\sigma$, the condition:
\par
\noindent a)$\qquad t_1\in T_{P_1}(S)$
\par
\noindent does not imply the condition :
\par
\noindent b)$\qquad t_2\in T_{P_2}(S)$ .
\par
\medskip
Fix the basis $\{F_1,...,F_{\sigma+1}\}$ for $(I_C)_\sigma$,
so we can write any $f\in (I_C)_\sigma$ as:
\par
$f= \lambda_1F_1+\lambda_2F_2+...+\lambda_{\sigma +1}F_{\sigma +1} =
{\rm det} \pmatrix {&&M&\cr
\lambda_1&\lambda_2&...&\lambda_{\sigma +1}}$
\par
\noindent and note that condition a) determines a
$\sigma$-dimensional vector space ${\bf V}_1\subseteq (I_C)_\sigma$ .
\par
Let :
$$\qquad M_k = \pmatrix{\delta^k_{11}&\delta^k_{12}&\quad ...\quad
 &\delta^k_{1\sigma +1}\cr
\delta^k_{21}&\delta^k_{22}&\quad ...\quad &\delta^k_{2\sigma +1}
\cr \vdots & \vdots&&\vdots \cr
\delta^k_{\sigma 1}&\delta^k_{\sigma 2}&\quad ...\quad
&\delta^k_{\sigma \sigma +1}}; $$
It is well-known
and easy to check that the rows of $M_0$ give elements
$f=\delta^0_{i1}F_1+...+\delta^0_{i\sigma +1}F_{\sigma +1} \in {\bf V}_1$,
since $f = {\rm det} \pmatrix {&&M&\cr
\delta^0_{i1}&\delta^0_{i2}&...&\delta^0_{i\sigma +1}}$,  and this
matrix has rank
$\sigma - 2$ at $P_1$ (rank $M_0={\rm rank}M(P_1) = \sigma -1$,
 since $P_1 \in C$ and $C$ is smooth).
 \par
In an analogous way, b) determines a $\sigma$-dimensional subspace
${\bf V}_2 \subseteq (I_C)_\sigma$,
which contains the vectors given by the rows,
$\underline \delta^3_{i}$, of $M_3$.
\par
There exists an element $\underline a =
(a_1,...,a_{\sigma +1}) \in k^{\sigma +1}$ such that the
set
${\cal B} = \{\underline \delta^3_{1},...,
\underline \delta^3_{\sigma},\underline a\}$
is a (redundant) system of generators of ${\bf V}_2$.
\par
\medskip
If a) $\Rightarrow $ b), then ${\bf V}_1 \subseteq {\bf V}_2$, and
we have that all vectors
$\underline \delta^0_{i}$ are linear combinations of elements of
${\cal B}$, i.e. $\forall i=1,...,\sigma$ there exists
$(\alpha_1^i,...,\alpha_{\sigma}^i,\alpha ^i)\in k^{\sigma +1}$ such that:
\par
\noindent$(\ddagger ) \qquad \qquad
\underline \delta^0_{i} = \sum_{k=1}^\sigma \alpha^i_k
\underline \delta^3_{k} \quad+\quad \alpha^i\underline a.$
\par
Now consider the line $L=\{w_1=w_2=0\}$: we have that
$M\vert _L = w_0M_0 + w_3M_3$.  Since ${\rm rank}M_3 = \sigma -1$,
 there will be two independent linear combinations of the columns of
 $M_3$ which give the zero-vector:
$$\beta_1 \pmatrix {\delta^3_{11}\cr\vdots \cr
 \delta^3_{\sigma 1}} + ...+
\beta_{\sigma +1} \pmatrix {\delta^3_{1\sigma +1}\cr\vdots \cr
\delta^3_{\sigma \sigma+1}} =
\pmatrix {\sum_{j=1}^{\sigma +1} \beta_j \delta^3_{1j}\cr\vdots \cr
 \sum_{j=1}^{\sigma +1} \beta_j \delta^3_{\sigma j}} = \pmatrix {0\cr\vdots
 \cr 0};$$
$$\gamma_1 \pmatrix {\delta^3_{11}\cr\vdots \cr
\delta^3_{\sigma 1}} + ...+
\gamma_{\sigma +1} \pmatrix {\delta^3_{1\sigma +1}\cr\vdots \cr
\delta^3_{\sigma \sigma+1}} =
\pmatrix {\sum_{j=1}^{\sigma +1} \gamma_j \delta^3_{1j}\cr\vdots \cr
 \sum_{j=1}^{\sigma +1} \gamma_j \delta^3_{\sigma j}} = \pmatrix {0\cr\vdots
 \cr 0}.$$
If we use $(\beta_1,...,\beta_{\sigma +1})$ as coefficients
for a linear combination of
the columns of $M\vert_L$, we get:
$$
\beta_1 \pmatrix {\delta^0_{11}w_0+ \delta^3_{11}w_3\cr\vdots \cr
\delta^0_{\sigma 1}w_0+\delta^3_{\sigma 1}w_3} + ...+
\beta_{\sigma +1} \pmatrix {\delta^0_{1\sigma +1}w_0 +
\delta^3_{1\sigma +1}w_3\cr\vdots \cr \delta^0_{\sigma \sigma +1}w_0 +
\delta^3_{\sigma \sigma+1}w_3} =
w_0 \pmatrix {\sum_{j=1}^{\sigma +1} \beta_j \delta^0_{1j}
\cr\vdots \cr
 \sum_{j=1}^{\sigma +1} \beta_j \delta^0_{\sigma j}}. $$
By $(\ddagger )$, we have:
$\delta^0_{ij} = \sum_{k=1}^\sigma \alpha^i_k
\delta^3_{kj} \ +\ \alpha^i a_j$, hence:
$$
w_0 \pmatrix {\sum_{j=1}^{\sigma +1}
\sum^\sigma_{k=1}\beta_j(\alpha_k^1\delta^3_{kj} + \alpha^1a_j)
\cr\vdots \cr
 \sum_{j=1}^{\sigma +1}
\sum^\sigma_{k=1}\beta_j (\alpha^\sigma_k\delta^3_{kj} + \alpha^\sigma a_j)}
=
w_0 \pmatrix {\sum^\sigma_{k=1}\alpha^1_k
(\sum_{j=1}^{\sigma +1}(\beta_j \delta^3_{kj})+
\alpha^1\sum_{j=1}^{\sigma +1}\beta_ja_j
\cr\vdots \cr
\sum^\sigma_{k=1}\alpha^\sigma_k(\sum_{j=1}^{\sigma +1} \beta_j \delta^3_{kj})
+  \alpha^\sigma\sum_{j=1}^{\sigma +1} \beta_ja_j}$$
which, since $\sum_{j=1}^{\sigma +1} \beta_j \delta^3_{kj}=0$, yelds:
$$
(\sum_{j=1}^{\sigma +1}\beta_ja_j )w_0\pmatrix {\alpha^1
 \cr\vdots \cr
\alpha^\sigma}.$$
In an analogous way one gets:
$$
w_0 \pmatrix {\sum_{j=1}^{\sigma +1}
\sum^\sigma_{k=1}(\gamma_j \alpha_k^1\delta^3_{kj} + \alpha^1a_j\gamma_j)
\cr\vdots \cr
 \sum_{j=1}^{\sigma +1}
\sum^\sigma_{k=1}(\gamma_j \alpha^\sigma_k\delta^3_{kj} +
\alpha^\sigma a_j\gamma_j)}=
(\sum^{\sigma +1}_{j=1}\gamma_ja_j )w_0\pmatrix {\alpha^1
 \cr\vdots \cr
\alpha^\sigma}.$$
Since these two columns are multiples of each other, by
elementary columns operations we can replace a column in $M\vert _L$,
say the first one, by a column of zeros.
\par
Hence either the rank of $M$ drops at each point of $L$, i.e. $L\subseteq
C$, or
 all the minors of $M\vert_L$ are zero except  $F_1\vert _L$, i.e.
 $C\cap L$ consists of $\sigma$ points (the zeros of $F_1\vert _L$).
 \par
Since both cases cannot occur, we are done.
\par\medskip
In the case when $L_1\cap L_2 \neq \emptyset$, they are contained in
 a plane $\Gamma$, say $\Gamma = \{w_3 =0\}$. Then
the procedure is exactly the same, simply choose
on $\Gamma$: $P_1=(1:0:0)$, $P_2=(0:0:1)$, $L_1=\{w_2=0\}$,
$L_2=\{w_0=0\}$, $L=\{w_1=0\}$.
\par
\bigskip
\noindent {\bf Case 2)}.  This case is quite similar to the previous one.
Let $P_i=\pi(Q_i)$, so $P_1\in C$, $P_2\notin C$, and let
$t_1\in T_{P_1}(\PP 3) - T_{P_1}(C)$ be as in the previous case.
 \par
As before, we have to check that,
for $S=\{f=0\}$, $f\in (I_C)_\sigma $, the condition
\par
\noindent a) $\quad t_1 \in T_{P_1}(S)$
\par
\noindent does not imply the condition:
\par
\noindent b) $\quad P_2 \in S$.
\par
Choose for $P_1,P_2,L_1$ the same coordinates as before.
Since b) gives a subspace
of $(I_C)_\sigma$ defined by
$\{\underline \delta^3_1,...,\underline \delta^3_{\sigma+1}\}$ (this
time they are independent since $P_2\notin C$, i.e. rk$M_3=\sigma$), if
a) $\Rightarrow$ b), we will again have that all the rows in
$M_0$ are linear combination of the rows in $M_3$, and we can conclude as
before (here we just have $\underline a = \underline 0$).
\par
\rightline {\qed}
\bigskip
It is possible to extend the result above to the case of varieties of
codimension 2 which do not contain lines.
\bigskip
\noindent {\bf Proposition 2.4:}  {\it Let $V\subseteq \PP n$
be a smooth a.C.M. variety of codim 2 and let
$I_V$ be generated by the maximal minors of a matrix of linear forms.
Then} $\vert \sigma H-E\vert$ {\it is very ample on} $X_V$ {\it if
V does not  contain any line and it has no} $\sigma${\it -secant line.}
\par
\medskip
The proof of Prop. 2.4 follows the same  lines as that of Thm. 2.3;
the hypothesis
that $V$ does not contain lines is needed to exclude the possibility that
the rank of $M$ drops at every point of $L$ (see Case 1 in the proof of 2.3).
\par
\rightline {\qed}
\bigskip
Now, the next problems that are quite natural to consider for these
varieties are:
\par
\noindent 1) when does $V$, defined by a generic matrix of
linear forms, contain lines?
\par
\noindent 2) when does such a $V$ have $\sigma$-secants?
\par
\medskip
Note that here "generic" matrix means that the coefficients of its entries
are generic in $k^{\sigma(\sigma +1)(n+1)}$.
\par
\medskip
The answer to these questions is given by the next proposition
and it was probably classically known. We give a proof here,
 for lack of reference (notice that we consider only the case
$\sigma \geq n+1$, since for $\sigma \leq n$ $V$ trivially has
$\sigma$-secants).
\par
\bigskip
\noindent {\bf Proposition 2.5:}  {\it Let M be a generic}
$\sigma \times (\sigma +1)${\it -matrix,} $\sigma \geq n+1${\it ,
whose entries are linear forms in}
$k[w_0,w_1,w_2,...,w_n]$, $n\geq 3$.
{\it Then M defines an irreducible
a.C.M. variety V in} $\PP n$ {\it and }
\par
\noindent {\it i)  V contains a line if and only if} $\sigma \leq 2n-5$.
\par
\noindent {\it ii)  V has } $\sigma${\it -secants if and only if}
$\sigma \leq 2n-2$.
\par
\medskip
\noindent {\it Proof}: For the fact that
a generic matrix defines an
irreducible 2-codimensional variety see
e.g. \ref S (recall, see the introduction,
that $V$ will be smooth if and ony if $n\leq \sigma$).
\par
As before, let
$$\qquad M = \pmatrix{L_{11}&L_{12}&\quad ...\quad &L_{1\sigma +1}\cr
L_{21}&L_{22}&\quad ...\quad &L_{2\sigma +1}
\cr \vdots & \vdots&&\vdots \cr
L_{\sigma 1}&L_{\sigma 2}&\quad ...\quad &L_{\sigma \sigma +1}}\qquad
{\rm where}\qquad L_{ij}=\sum_{k=0} ^n
\delta_{ij}^kw_k .$$
  For each $\underline{z}=[z_1:...:z_{\sigma}]$, let us consider
the $\sigma + 1$ hyperplanes:
\par
\noindent ($\dagger$) $\qquad  \qquad
\{z_1L_{1j} + z_2L_{2j} + ... + z_\sigma L_{\sigma j} = 0\},
\qquad j=1,...,\sigma +1$.
\par
Then $V$ can be viewed as the locus of meets of the $\sigma + 1$
hyperplanes above, as $\underline{z}$ varies in $\PP {\sigma -1}$ (this is
classically known as the {\it projective generation} of $V$).
\par
\medskip
\noindent Proof of $i)$:
\par
For a given $\underline{z} \in \PP {\sigma -1}$,
let us denote by $T_{\underline{z}} \subset V$ the linear space of the
solutions of the system (in the $w_k$'s) defined by ($\dagger$), whose
matrix of coefficients $Z$ is an $(n+1)\times (\sigma +1)$-matrix of
linear forms in $k[z_1,...,z_\sigma ]$:
$$
Z_{k j} = \sum _{i=1} ^\sigma \delta _{i j} ^k z_i.
$$
We will have that $T_{\underline{z}}$ contains a line if and only if
the rank of $Z$ is $\leq n-1$. $Z$ is a generic matrix since $M$ is (the
coefficients that show up are the same), hence the codimension of the
locus $\Lambda $ where the rank of $Z$ is $\leq n-1$ is $(\sigma +1 -n
+1)(n+1-n+1) =
 2\sigma -2n +4$. Thus $\Lambda \subset \PP {\sigma -1}$
 is not empty if and only if   $2\sigma -2n +4 \leq \sigma - 1$, i.e.
whenever $\sigma \leq 2n-5$.
\par
\medskip
Now let $L$ be a line in $V$; if we show that $L$ is contained in some
$T_{\underline{z}}$ we are done, but this is obvious
since $L \subseteq V$, so the rank of $M\vert _L$ drops, i.e. its rows
are linearly dependent, therefore there is a solution
$\underline {z}$ for the system $(\dagger )$, hence
$L \subseteq T_{\underline{z}}$.
\par
\medskip
\noindent Proof of $ii)$:
\par
Every family of hyperplanes given by
$z_1L_{1j}+z_2L_{2j}+...+z_{\sigma}L_{\sigma j} = 0$ has a "center",
given by $\{L_{1j}=L_{2j}=...=L_{\sigma j} = 0\}$. Via a linear
combination of the columns of $M$,
one can consider
$$
T_{\underline y} =
\{\sum_{j=1}^{\sigma +1} y_jL_{1j}=\sum_{j=1}^{\sigma +1} y_jL_{2j}=...
=\sum_{j=1}^{\sigma +1 } y_jL_{\sigma j} = 0\}
$$
\noindent
as one of the "centers" of ($\dagger$), for every choice of
$\underline y\in \PP \sigma$.
\par
If we consider the linear system defined by $T_{\underline{y}}$
in the $w_k$'s, the entries of its
$(\sigma \times(n+1))$-matrix of coefficients $N$ are given by
the forms (in $k[y_1,...,y_{\sigma +1}]$):
$$
N_{ik} = \sum_{j=1}^{\sigma +1} \delta_{ij}^k y_j
$$
Now, for a given $\underline{y}$, if the matrix $N(\underline y) $
has rank $n+1$, then $T_{\underline y}$ is empty; if rk$N(\underline y)
=n$, then
$T_{\underline y}$ is a point; while $T_{\underline y}$ is a line for
rk$N(\underline y) = n-1$.
\par
In the last case, the line $T_{\underline {y}}$ intersects $V$
(see e.g. \ref R)
in a subvariety of $T_{\underline y}\cong \PP 1$ defined by a
$(\sigma\times \sigma)$-matrix of linear forms,
i.e. in $\sigma$ points (with
multiplicities, i.e. in a divisor of degree $\sigma$ in $\PP 1$).
\par
Since $M$ was generic, so is $N$ (the coefficients of the forms are the
same), hence the scheme $\Gamma \subseteq \PP \sigma$, defined as the
rank $n-1$ locus of $N$, has codimension
$ (\sigma -(n-1))\times(n+1-(n-1)) = 2(\sigma -n+1) = 2(\sigma -n+1)$ in
$\PP \sigma$, so
it is not empty for $\sigma \leq 2n-2$, as we wanted.
\par
\medskip
\noindent   We are left with checking that when $\sigma \geq 2n-1$,
there are no $\sigma$-secants to $V$.
\par
We will be done if we show that, when $V$ possesses a $\sigma$-secant $L$,
$L$  can be viewed as before, i.e. as the
"center" of a family of hyperplanes
which gives the projective generation of $V$.
\par
Let $L = \{w_2=w_3=...=w_n=0\}$, then on $L$ we have that
$M\vert _L$ defines the
homogeneous ideal, in $R\prim =k[w_0,w_1]$, of $\sigma$ points, i.e. that its
 maximal
minors $F_1,...,F_{\sigma +1}$ (where $F_j$ is the minor with sign
obtained by erasing the $j^{th}$ column)
can be viewed as multiples of one of them
(which is not identically zero), say $F_1$.
\par
This implies that the columns of the matrix $M\vert _L$,
are linearly dependent, hence there is a linear combination of them,
say with coefficients $y_1,...,y_{\sigma +1}$,
which gives a zero column.
\par
Then the same linear combination on $M$ gives a form in $k[w_2,w_3,...,w_n]$,
i.e. we have that:
$$\{\sum_{j=1}^{\sigma +1 } y_jL_{1j}=\sum_{j=1}^{\sigma +1} y_jL_{2j}=...
=\sum_{j=1}^{\sigma +1} y_jL_{\sigma j} = 0\}$$
\noindent
has solution on all of $L= \{w_2=w_3=...=w_n=0\}$. But the matrix of this
linear system is what we have called $N(\underline{y})$, so it must have
rank $\leq n-1$, which is possible only for $\sigma \leq 2n-2$,
as we have observed before.
\par
\rightline {\qed}
\bigskip
\medskip
Let us give an example of what happens for the first values of
$\sigma $,  when $V=C$ is a curve:
\par
- when $\sigma = 2$, $C=C^3_0$ is a rational normal curve and, for
every $\underline {y} \subseteq \PP 2$, $T_{\underline y}$ is a secant line
of $C$.
\par
- when $\sigma =3$, $C=C^6_3$ is a sestic curve of genus three, and $\Gamma
\cong C$
gives a family of trisecants of $C$.
\par
- when $\sigma =4$, $\Gamma \subset \PP 4$ is a set of 20 points,
corresponding to 20 4-secants of $C=C^{10}_{11}$.
\par
\bigskip
>From the above results, it follows the analog of the Conjecture
for a "generic" smooth $V$ of codimension 2:
\par
\bigskip
\Coroll {2.6}:
{\it Let $M$ be a generic $\sigma \times (\sigma +1)$-matrix
of linear forms in $k[w_0,...,w_n]$, with $1<n<6$, and let
$V\subseteq \PP n$
be the smooth a.C.M. variety of codim $2$ defined by
$I_V$, the ideal generated by the maximal minors of $M$.
Then} $\vert \sigma H-E\vert$ {\it is very ample on} $X_V$ {\it if and
only if} $\sigma > 2n-2$.
\par
\vskip 1cm
\noindent {\bf 3. Generic curves of minimal genus.}
\medskip
In this section we will consider projectively normal curves
$C \subseteq \PP 3$ such that they are
"generic for their degree" in the sense that if $\deg C = s$, the
ideal sheaf of $C$ has the same kind of resolution (in the sense
that it has the same graded Betti numbers) as a generic set
of $s$ points in $\PP 2$. Which are the Betti numbers of this
"generic resolution" (in the sense of \ref {L.2}) will be specified
below (see comments before Prop. 3.3).
\par
We will say that such a $C$ is {\it Betti minimal},
B-minimal for short, and note that its Hilbert-Burch matrix
has entries of minimal degree. As we will see below, such curves have
minimal genus among projectively normal curve of their degree.
\par
Since a generic hyperplane section of $C$ has the same
graded Betti numbers as $C$, by cutting with a generic plane
in order to understand the Hilbert-Burch
matrix of a B-minimal $C$, we can reduce ourselves
to look at a (reduced) set of
points $Z\in \PP 2$. The requirement on the resolution
of $Z$ implies that the degree matrix of the
ideal $I_Z$ (which is the same as the degree matrix of $I_C$) is made
by 1's and 2's and the Hilbert function of $Z$ is maximal, i.e.
$H(Z,t) = min\{{t+2 \choose 2},s\}$ (see \ref {C-G-O},\ref{Gi-Lo}).
\par
The Hilbert function of $Z$ is related to the
genus of $C$ by the following well-known fact (which we prove here
for lack of reference):
\par
\bigskip
\noindent {\bf Proposition 3.1:} {\it Let $C\subseteq \PP 3$ be
projectively normal, and let Z be
a generic plane section of C. Let} $\deg C = s$; {\it
then the genus of C is:}
$g = \sum _{t\geq 1} s - H(Z,t)$.
\par
\medskip
\noindent {\it Proof:} We have $g = h^2(\PP 3,{\cal I}_C)$, and
$s - H(Z,t) = h^1(\PP 2,{\cal I}_Z(t))$.
>From the resolution of $C$:
\par
$$0 \rightarrow \oplus _{j=1}^u{\cal O}_{\PP 3}(-n_j) \rightarrow
\oplus _{i=1}^v{\cal O}_{\PP 3}(-d_i) \rightarrow {\cal I}_C \rightarrow 0,$$
\noindent we get that
\par
\medskip
\noindent $h^2({\cal I}_C) =
\sum_{j=1}^u h^3(\PP 3,{\cal O}_{\PP 3}(-n_j)) -
\sum_{i=1}^v h^3(\PP 3,{\cal O}_{\PP 3}(-d_i)) = $
\par
$ = \sum_{j=1}^u h^0(\PP 3,{\cal O}_{\PP 3}(n_j-4))-
\sum_{i=1}^v h^0(\PP 3,{\cal O}_{\PP 3}(d_i-4)) =
\sum_{i=1}^u {n_j-1 \choose 3} - \sum_{i=1}^v {d_i-1 \choose 3}$.
\medskip
The s points $Z\subseteq \PP 2$ have the same graded Betti numbers as $C$.
Hence, for their ideal sheaf
${\cal I}_Z \subseteq {\cal O}_{\PP 2}$, we have:
\par
\medskip
\noindent $h^1(\PP 2,{\cal I}_Z(t)) =
\sum_{j=1}^u h^2(\PP 2,{\cal O}_{\PP 2}(t-n_j)) -
\sum_{i=1}^v h^2(\PP 2,{\cal O}_{\PP 2}(t-d_i)) = $
\par
$\qquad  = \sum_{j=1}^u h^0(\PP 2,{\cal O}_{\PP 2}(n_j-t-3))-
\sum_{i=1}^v h^0(\PP 2,{\cal O}_{\PP 2}(d_i-t-3)) = $
\par
\centerline
{$ \sum_{j=1}^u {n_j-t-1 \choose 2} - \sum_{i=1}^v {d_i-t-1 \choose 2}.$ }
\par
\medskip
Let us consider now:
\par
\medskip
\noindent $\sum _{t\geq 1} h^1(\PP 2,{\cal I}_Z(t)) =
\sum _{t\geq 1}(\sum_{j=1}^u {n_j-t-1 \choose 2} -
\sum_{i=1}^v {d_i-t-1 \choose 2}).$
\par
\medskip
Since $h^1(\PP 2,{\cal I}_Z(t)) = 0$, for $t\geq \sigma -1 = max_j\{n_j\}-2
> max_i\{d_i\}-2$,
if we put, in any sum over t, $\alpha =n_j-t-1$, and
$\beta = d_i-t-1$, then we eventually get
\par
\medskip
$\noindent \sum _{t\geq 1} h^1(\PP 2,{\cal I}_Z(t))=
\sum_{j=1}^u (\sum _{\alpha = 1}^{n_j-2} {\alpha \choose 2}) -
\sum_{i=1}^v (\sum _{\beta = 1}^{d_i-2} {\beta \choose 2}) =$
\par
\centerline{$
\sum_{j=1}^u {n_j-1 \choose 3} -
\sum_{i=1}^v {d_i-1 \choose 3} = g(C)$}
\medskip
\noindent and we are done.
\par
\rightline {\qed}
\bigskip
As an immediate consequence of the proposition above and the fact
that B-minimal curves have plane sections with maximal Hilbert
function, we get:
\par
\medskip
\noindent {\bf Corollary 3.2:} {\it Any B-minimal curve has minimal
genus among the projectively normal curves of its degree.}
\par
\medskip
Now let us describe more precisely the degree matrix of B-minimal curves;
let $s=\deg C = \deg Z$,
$ d = min \{t\vert s\leq {t+2 \choose 2}\}$ and $s = {d+1 \choose 2} + k$,
with $1 \leq k \leq d+1$ (note that we have $\sigma = \sigma (I_C)=d+1$).
\par
The ideal $I_C$ will be generated by $d-k+1$ forms of degree $d$ and $h$ forms
of degree $d+1$, where $h$ is either $0$ or $2k-d$,
according to whether $d\geq 2k$ or not.
\par
Since given any B-minimal $C$ the integers $d,k$ are
determined by its degree $s$, we will also speak of
"B-minimal curve of type $d,k$".
\par
Denote by $F_1,...,F_{d-k+1}$ the generators of degree $d$ and by
$G_1,...,G_{2k-d}$ those of degree $d+1$ (when present).
\par
We have that the Hilbert-Burch matrix is a $\rho\times (\rho +1)$- matrix,
where :
$$ \rho = \left\{ \matrix { k &\ &{\rm if}\  d\leq 2k \cr
d-k &\ &{\rm if}\  d\geq 2k }\right .$$
In the case $d\leq 2k$, the matrix $M$ is :
$$M = \pmatrix{L_{11}&L_{12}&...&L_{12k-d}
&Q_{11}&...&Q_{1d-k+1}\cr
L_{21}&L_{22}&... &L_{2 2k-d} &Q_{21}&...&Q_{2d-k+1}
\cr \vdots& \vdots&&\vdots &\vdots&&\vdots\cr
L_{k1}&L_{k2}&\quad ...\quad &L_{k2k-d}&Q_{k1}&...&Q_{kd-k+1} }.$$
\noindent where the $L_{ij}$'s are linear forms and the $Q_{il}$'s
are forms of degree 2.  The $F_j$'s are the algebraic minors
obtained by deleting the $j^{th}$ column, while the $G_l$'s
are obtained by deleting the $l^{th}$ column.
\par
\medskip
In this case the resolution of $I_C$ is:
$$0 \rightarrow \oplus ^{k}{\cal O}_{\PP 3}(-d-2) \rightarrow
(\oplus ^{2k-d}{\cal O}_{\PP 3}(-d-1))
\oplus (\oplus ^{d-k+1}{\cal O}_{\PP 3}(-d)) \rightarrow {\cal I}_C \rightarrow
0.$$
\medskip
In the other case ($d\geq 2k$) the matrix is:
$$ M = \pmatrix{Q_{11}&Q_{12}&... &Q_{1d-k+1}\cr
\vdots& \vdots&&\vdots \cr
Q_{k1}&Q_{k2}& ... &Q_{kd-k+1}\cr
L_{11}&L_{12}&  ...  &L_{1d-k+1}
\cr \vdots& \vdots&&\vdots \cr
L_{d-2k 1}&L_{d-2k 2}& ... &L_{d-2k d-k+1}}. $$
\medskip
Hence the resolution of ${\cal I}_C$ is:
$$0 \rightarrow (\oplus ^{k}{\cal O}_{\PP 3}(-d-2))
\oplus (\oplus ^{d-2k}{\cal O}_{\PP 3}(-d-1))
\rightarrow
\oplus ^{d-k+1}{\cal O}_{\PP 3}(-d) \rightarrow {\cal I}_C \rightarrow 0.$$
\medskip
In this case we have that the $F_j$ are the maximal algebraic
minors of $M$, and the ideal $I_C$ is generated
in degree $d$.
\par
\medskip
>From what we have observed, it easily follows that our conjecture is
trivially true
for B-minimal curves of type $d$,$k$ with $d\geq 2k$:
\par
\bigskip
\noindent {\bf Proposition 3.3:} {\it Let C be B-minimal of type d,k with}
$d\geq 2k${\it .  Then}
$\vert \sigma H -E\vert $ {\it is very ample on} $X_C$.
\par
\medskip
\noindent {\it Proof:} The Proposition follows immediately from
\ref {Co}, Theorem 1, quoted in the introduction, since $I_C$ is
generated in degree $d = \sigma -1 $ (notice that for this reason
$C$ cannot possess $\sigma$-secants).
\par
\rightline {\qed}
\bigskip
Let us notice that, by what we have seen in the previous section,
for any $C$ of type (*), the linear system
$\vert \sigma H -E\vert $ always defines a
morphism $\phi _\sigma :X_C \rightarrow \PP N$, where
$N+1 = h^0(X_C,{\cal O}_{X_C}(\sigma H -E))= \dim (I_C)_\sigma$,
whose image
is a rational threefold (by Prop.1.1, $\phi _\sigma $ is generically
1:1). We have:
\par
\bigskip
\noindent {\bf Theorem 3.4:} {\it Let C}$\subseteq \PP3$
{\it be a projectively normal curve
which is B-minimal of type d,k. Then there exist two matrices B,X of
linear forms (in the coordinates ring S
of} $\PP N${\it ) such that: X is a} $4\times (d-k+1)${\it -matrix, B is
a} $k\times 4${\it -matrix, and the homogeneous ideal I of}
$Y_{C,\sigma} = \phi _\sigma (X_C)$
{\it in S is generated by the}
$(2\times2)${\it -minors of X, the entries of} $B\cdot X$ {\it and the }
$(4\times 4)${\it -minors of B.  Moreover, $Y_{C,\sigma} $ is a.C.M.}
\par
\medskip
\noindent {\it Sketch of proof:}
The proof is very similar to the one given in $\S 4$ of \ref {Gi-Lo},
for points in $\bf P^2$, which essentially relies upon the structure
of the matrix $M$ and its genericity, hence it easily extends to this
case.  We refer to that paper for full details, here we simply
explain the main steps of the procedure:
\par
1) Construction of a homogeneous ideal $I$ as described in the statement,
such that $Y_{C,\sigma } \subseteq Z(I).$
\par
2) Prove that $I$ is prime and perfect.
\par
3) Show that $Y_{C,\sigma }=Z(I).$
\par
\noindent {\bf  Step 1} - The idea for the construction of $I$ is the
fact that its elements can be obtained from syzygies of $I_C$ whose
components belong to $I_C$ itself; so we have to show that we can find
the appropriate syzygies which will give us a set of generators for $I.$
\par
First of all, recall that we denoted by $S=k[\underline{x},\underline{y}]$
the coordinate ring of $\bf P^N$, where $\underline{x} = (x_{hj})$,
$h=0,...,4$; $j=1,...,d-k+1$,
and by  $R=k[\underline{w}]$ the coordinate ring of $\bf P^3$.
Let $R\prim$ be the graded subring of $R$ defined by
$R\prim =\oplus _k R\prim _k,$ where $R\prim _k=R_{k(d+1)}.$
Now define $\psi : k[\underline{x},\underline{y}]
\rightarrow k[\underline{w}]$ by $\psi (x_{hj})=(w_hF_j),$
 $\psi (y_l)=G_l$, which is graded (i.e. homogeneous of degree $0$).
Let $N\prim = 3d-3k+\rho +3.$  When $d\leq 2k,$ we have that
$N\prim=N$ and $I_Y = Ker \psi.$  When $d>2k,$ we can derive, from
the matrix $M,$ $d-2k$ linearly independent linear forms
$H_1,...,H_{d-2k}\in Ker \psi$ so that $I_Y$ may be identified with
${ker\psi \over (H_1,...H_{d-2k})}.$
\par
The first obvious generators of $I$ will be all the elements
$x_{h\prim j\prim}x_{hj}-x_{h\prim j}x_{hj\prim}$, which are in $Ker \psi$
(they correspond to the trivial syzygies
$w_{h\prim}F_{j\prim}w_hF_j-w_{h\prim}F_jw_hF_{j\prim}$),
and which can be viewed as the $(2\times 2)$-minors of the matrix:
$$X=\pmatrix{x_{01}&x_{02}&...&x_{0d-k+1}\cr
x_{11}&x_{12}&...&x_{1d-k+1}\cr
&&...&\cr x_{31}&x_{32}&...&x_{3d-k+1}}.$$
\par
Now assume $d\leq 2k$ and consider the matrix $M_u$, $u=1,...,k$, obtained by
 repeating
the $u^{th}$ row in $M$:
$$M_u = \pmatrix{&&&M&&&\cr
 L_{u1}&L_{u2}&...&L_{u2k-d}&Q_{u1}&...&Q_{ud-k+1}};$$
\noindent
let $L_{ul} = \sum_{i=0}^{3}\delta^{ul}_iw_i$, and
$Q_{uj} = \sum _{i,h=0}^3\beta^{uij}_hw_hw_i$; by
expanding $\det M_u $ with respect to the last row we get:
$$0 = \sum_{l=1}^{2k-d}L_{ul}G_l + \sum_{j=1}^{d-k+1} Q_{uj}F_j =
\sum_{l=1}^{2k-d}\left(\sum _{i=0}^3\delta^{ul}_iw_i\right)G_l +
\sum_{j=1}^{d-k+1}\left( \sum _{i,h=0}^3\beta^{uij}_hw_hw_i\right)F_j;$$
\noindent
and so, by multiplying by $F_\nu$, we get $\forall \nu = 1,...,d-k+1$:
$$0 = F_\nu \det M_u =
\sum _{i=0}^3\left(\sum_{l=1}^{2k-d}\delta^{ul}_i\psi (y_l)\right)\psi
 (x_{i\nu}) +
\sum _{i}^3\left(\sum_{j=1}^{d-k+1} \sum _{h=0}^3\beta^{uij}_h\psi (x_{hj})
\right)\psi (x_{i\nu})$$
which can be written as:
$$0= \psi \left(\sum _{i}^3x_{i\nu}B_{ui}\right),$$
where
$$B_{ui} = \sum_{l=1}^{2k-d}\delta^{ul}_iy_l +
\sum_{j=1}^{d-k+1} \sum _{h=0}^3\beta^{uij}_h x_{hj}.$$
The degree 2 forms $\sum _{i}^3x_{i\nu}B_{uj}$ in $\ker \psi$ will
be added as generators of $I$ and
can be viewed
as the entries of the product matrix $B\cdot X$, where $B$ is the
$(4\times k)$-matrix of linear forms : $B = \left(B_{ui}\right)_{ui}$.
\par
Furthermore, denote $C_{ui}=\psi (B_{ui});$ then  the image under
$\psi$ of the maximal minors of $B$ are precisely the maximal minors
of the matrix $C=(C_{ui})_{ui}.$  Now, for each $u=1,...,k$
and $i=0,...,3$, consider the matrix
$$D_{ui}=\pmatrix{&&&M&&&\cr
\delta ^{ul}_i&...&\delta^{u2k-d}_i&\sum _{h=0}^3\beta ^{uh1}_iw_h&...&
\sum _{h=0}^3\beta ^{uhd-k+1}_iw_h}.$$
Then, it turns out that $\det D_{ui}=C_{ui}$ and that, $\forall
(a_0:...:a_3) \in \PP 3$,
$$0=\det M_u(a_0:...:a_3)=\sum _{i=0}^3a_i\det D_{ui}(a_0:...:a_3)=
\sum _{i=0}^3a_iC_{ui}(a_0:...:a_3).$$
Therefore the rank of $C$ is not maximal, and so all the maximal
minors of $B$ belong to $Ker \psi .$   We add these forms of degree
$4$ to the generators of $I.$
\medskip
The case $d>2k$ is alike, but the matrix $M_u,$ obtained by repeting the
$u^{th}$ row of $M$, in this case looks like
$$M_u = \pmatrix{&M&\cr
Q_{u1}&...&Q_{ud-k+1}}.$$
\par
Note that,  in this case, there are no $G_l$, hence no $y_l$ will show up
in the resulting matrix $B.$
\par
Also, in this case, $I$ will be defined as the ideal generated
 by the residue classes
of the $2 \times 2$ minors of $X$, of the entries of $B\cdot X$,
and of the $4 \times 4$ minors of $B$, modulo $(H_1,...,H_{d-2k})$.
\par
\medskip
\noindent {\bf Step 2} - A genericity argument and a Theorem by Huneke
(\ref {Hu, Thm. 60}, see also \ref {Gi-Lo, Thm. 4.1}) allow us to prove
that I is prime and perfect, just as in \ref {Gi-Lo, Thm. 4.2}.
\par
\medskip
\noindent {\bf Step 3} - After the necessary changes, the proof of
\ref {Gi-Lo , Prop. 4.4} works throughout.
\par
\rightline {\qed}
\bigskip
\noindent {\bf Remark 1.} When $k=d+1$, $M$ is a matrix of linear forms and
in this case only the matrix $B$
will appear in the construction of the ideal $I$ of $\phi _\sigma(C)$,
which will be determinantal (the situation will be the one described
in Corollary 2.2, with $k=\rho=\sigma$).
\par
\bigskip
\noindent {\bf Remark 2.} The threefold $Y_{C,\sigma}$ considered in 3.4 need
 not
 be smooth: there are many cases when $\phi _\sigma $ is not very ample (e.g.
see the Remark in the next section).
\par
\bigskip
There is another case in which it is not hard to compute the
ideal generation of the embedded threefold (this time via the
linear system $\vert (\sigma +1)H-E\vert$):
\par
\bigskip
\Prop {3.5} {\it Let } $C \subseteq \PP 3$ {\it be a projectively normal curve
whose Hilbert-Burch matrix is made of linear forms and let }
$\deg C = {d+1 \choose 2}$ {\it (hence } $\sigma (C) = d${\it ).
Then the homogeneous ideal of } $Y_{C,\sigma +1}\subseteq
\PP {3d+3}$ {\it is generated by the} $(2\times 2)$-{\it minors of a }
$4\times (d+1)${\it -matrix X of linear forms (in the coordinate ring of}
 $\PP {3d+3}${\it ).}
\par
\medskip
\noindent The proof works as in the previous theorem (where one could
consider $k=0$): this time only the matrix $X$ is to be considered (see
\ref {Ge-Gi} for a description of this construction in the case of points
in $\PP 2$).  In order to check that $Y_{C,\sigma +1}\subseteq
\PP N$, with $N= {3d+3}$, consider that $\deg C = {d+1 \choose 2}$, and
$g(C) = {d+1 \choose 2}({2d-5\over 3}) + 1$ (see Prop. 3.1), hence:
$$
H(C,d+1) = {d+1 \choose 2}(d+1) - g(C) + 1
$$
\noindent
implies
$$
N+1 = {d+1+3 \choose 3} - H(C,d+1) = {d+4 \choose 3} - {d+1 \choose 2}(d+1) +
 g(C) - 1 + 1 =  3d+4 .
$$
\rightline {\qed}
\bigskip
Notice that in this case also the resolution of $I_{Y_{C,\sigma +1}}$ can be
computed via a Lascoux complex (see \ref {La} and also \ref {PW}).
\par
\medskip
Finally, we can observe that it is not hard to extend the results of 3.4 and
3.5 to the case of 2-codimensional (smooth, irreducible, a.C.M.) subschemes
 $V\subseteq \PP n$,
 $n=4,5$.  Namely, let $V$ be as above, and also B-minimal (with an obvious
 extension of
the definition given for curves in $\PP 3$).  Let $\deg V = s = {d+1
\choose 2} + k$,
$0\leq k \leq d+1$, in analogy of what we did before.  Then, if $Y_{V,d+1}$ is
the image of $X_V$ via the linear system $\vert (d+1)H-E\vert $, we have:
\par
\bigskip
\Prop {3.6} {\it With the above notation and hypotheses,
we have that the homogeneous ideal of } $Y_{V,d+1} \subseteq \PP N$ {\it is
 generated by
the} $(2\times 2)$ {\it minors of a } $k\times (n+1)${\it -matrix of linear
 forms X, by
the }$(n+1)\times (n+1)$ {\it minors of a } $k\times (n+1)${\it -matrix B of
 linear forms
and by the entries of } $B\cdot X$. {\it Moreover, } $Y_{V,d+1}$ {\it is
a.C.M.}
\par
\medskip
The proof of 3.6 works as the ones of 3.4 and 3.5 .
\par
\rightline {\qed}
\vskip 1cm
\noindent {\bf 4.  Examples: Some Fano varieties
and their ideal generation.}
\bigskip
\noindent{\bf Example 1} Let $C = C^7_5 \subseteq \PP 3$ be a
projectively normal
curve of degree 7 and genus 5.  $C$ is B-minimal of type 3,1 and its
Hilbert-Burch matrix has degrees:
$$ \pmatrix{2&2&2\cr 1&1&1}. $$
Hence $I_C$ is generated by cubics, and in this case $\sigma = 4$, so
the linear system $\vert 4H-E\vert$ is very ample on $X_C$; notice that
$-K_{X_C} = 4H-E$, hence $X_C$ is a Fano threefold of index 1.
\par
The image of $X_C$ is a $Y_{C,4} \subseteq \PP {10}$, and
a curve-section of it can be viewed as the image of a curve $C\prim \subseteq
\PP 3$
which is the residual intersection,
with respect to $C$, of two quartic surfaces.
Such curve has degree 9 and its genus is given by the following
formula, where $a,b$ are the degrees of two surfaces whose intersection
is formed (scheme theoretically) by $C \cup C\prim$ (see e.g. \ref {P-S}) :
$$
g(C) - g(C\prim) = \left( {a+b \over 2} - 2\right)(\deg C - \deg C\prim)
$$
\noindent
hence $g(C\prim) = 9$ .
\par
So $Y_{C,4}$ has degree $2g(C\prim) - 2 = 16$ and sectional genus 9.
\par
By Theorem 3.4, the homogeneous ideal of
$Y_{C,4}$ will be generated by quadratic forms,
given by the $2\times 2$ minors of
a $(4\times 3)$-matrix $X$ of linear forms and by the entries of $B\cdot X$,
where $B$ is a $(1\times 4)$-matrix of linear forms.
\par
\bigskip
\noindent {\bf Remark.}
Notice that $C = C^7_5$ is the only B-minimal curve with $\sigma = 4$ for
which $\vert 4H-E\vert$ is very ample on $X_C$; the other B-minimal
curves with $\sigma = 4$ are in fact:
$C^8_7$, $C^9 _9$ and $C^{10}_{11}$, whose degree
matrices are, respectively:
$$
\pmatrix{1&2&2\cr 1&2&2}; \qquad  \pmatrix{1&1&1&2\cr 1&1&1&2\cr 1&1&1&2};
\qquad  \pmatrix{1&1&1&1&1\cr  1&1&1&1&1 \cr 1&1&1&1&1 \cr 1&1&1&1&1}.
$$
In the first case it is trivial that the two linear form in the first column of
the matrix define a 4-secant line; in the second it can be seen that
any smooth cubic surface containing the curve has six lines (one part of its
"double-six" among its 27 lines) which are 4-secants to the curve;
the last case has been considered in Proposition 2.5 (there are twenty
4-secant lines).
\par
\medskip
Notice also that in these cases $\vert 4H-E\vert$ is not even ample (any
of its multiples will have zero intersection with the 4-secants), hence
those $X_C$'s are not Fano.
\par
\bigskip
\noindent{\bf Example 2} Of course one can get a Fano threefold
also when $\sigma (C) < 4$; in particular
let us consider $C = C^6_3$: in this case $I_C$ is generated by the
$3\times 3$-minors of a $3\times 4$-matrix of linear forms and $\sigma (C)
= 3$,
hence we are in the case of Proposition 3.5: the Fano threefold
$Y_{C,4} \subseteq \PP {12}$  has degree 20 and sectional genus 11 (two
quartics through C have residual intersection in a $C^{10}_{11}$) and its
homogeneous ideal is generated by the $2\times 2$ minors of a
$(4\times 4)$-matrix of linear forms.
\par
\bigskip
\noindent{\bf Example 3} Working as in Example 1, it is
possible to find Fano varieties of dimension 4 and 5.
 In $\PP 4$, consider the generic determinantal surfaces
 of degree 11 and 12, whose
Hilbert-Burch matrices have, respectively, degrees:
$$
\pmatrix{2&2&2&2\cr 1&1&1&1 \cr 1&1&1&1}; \qquad  \pmatrix{2&2&2 \cr
2&2&2}.
$$
Those surfaces are generated in degree 4, hence the (anticanonical)
linear system $\vert 5H-E \vert$ is very ample on the blow-up $X_V$,
and the generation of the ideals of the embedded 4-folds is given
by Proposition 3.6.
\par
\medskip
Similarly, in $\PP 5$, we have that we can use
the two determinantal threefold (of degrees, respectively, 16 and 17)
whose Hilbert-Burch matrices have entries, respectively, of the following
 degrees:
$$
\pmatrix{2&2&2&2&2\cr 1&1&1&1&1 \cr 1&1&1&1&1 \cr 1&1&1&1&1}; \qquad
\pmatrix{2&2&2&2 \cr
2&2&2&2\cr 1&1&1&1}
$$
to get Fano 5-fold whose ideal generation is again
described by Prop. 3.6.
\vskip 1cm
\centerline{REFERENCES}
\bigskip
\noindent \ref {C-G-O}: C.Ciliberto, A.V.Geramita, F.Orecchia:
{\it Perfect Varieties with defining equations of high degree.}
 Boll.U.M.I. (7) {\bf 1-B} (1987), 663-647.
\par
\medskip
\noindent \ref{Co}:  M.Coppens: {\it Embeddings of blowing-ups.}
Seminari di Geometria 1991/93, Univ.di Bologna, Bologna (1994).
\par
\medskip
\noindent \ref{D-G}: E.Davis, A.V.Geramita: {\it Projective embeddings
of blow-ups of $\PP 2$.} Math. Ann. {\bf 279}, 435-448, (1988).
\par
\medskip
\noindent \ref{E-N} : J.A.Eagon,D.G.Northcott: {\it A note on the Hilbert
function of certain ideals which are defined by matrices.} Mathematika,
{\bf 9}, (1962), 118-126.
\par
\medskip
\noindent \ref {Ge-Gi}: A.V.Geramita, A.Gimigliano:
 {\it Generators for the defining ideal of certain rational sufaces .}
Duke Math. J. {\bf 62}, (1991), 61-83.
\par
\medskip
\noindent \ref{G-G-H}: A.V.Geramita, A.Gimigliano, B. Harbourne:
{\it Projectively Normal but Superaboundant Embeddings of  Rational Surfaces
 in Projective Spaces.}
 J. of Alg. {\bf 169}, (1994), 791-804.
\par
\medskip
\noindent \ref {Ge-Ma}: A.V.Geramita, P.Maroscia:
 {\it The ideal of forms vanishing at a finite set of points in }$\PP n$.
J. Algebra {\bf 90}, (1984), 528-555.
\par
\medskip
\noindent \ref {Gi.1}: A.Gimigliano : {\it On Veronesean Surfaces.}
Proc. Konin. Ned. Akad. van Wetenschappen, Ser. A, {\bf 92} (1989),
71-85.
\par
\medskip
\noindent \ref {Gi.2}: A.Gimigliano : {\it Our thin knowledge about
fat points.}
Queen's Papers in Pure and Applied Math. {\bf 83}, {\it The Curves Seminar
at Queen's, vol. VI.} (1989).
\par
\medskip
\noindent \ref {Gi.3}: A.Gimigliano : {\it  On Some Rational Surfaces.}
Seminari di Geometria {\bf '91-'93}, Univ. di Bologna (1994).
\par
\medskip
\noindent \ref {Gi-Lo}: A.Gimigliano, A.Lorenzini:
{\it On the ideal of Veronesean Surfaces.}
Can.J. Math. {\bf 45} (1993),758-777.
\par
\medskip
\noindent \ref {Hu}: C.Huneke: {\it The arithmetic perfection of
Buchsbaum-Eisenbud varieties and generic modules of projective dimension two.}
 Trans. Amer. Math. Soc. {\bf 265} (1981),211-233.
\par
\medskip
\noindent \ref {La}: A. Lascoux: {\it Syzygies des vari\'et\'es
determinantales}. Adv. in Math. {\bf 30} (1978), 202-237.
\par
\medskip
\noindent  \ref {L.1}: A.Lorenzini: {\it Betti numbers of perfect
homogeneous ideals.} J. of Pure and Applied Alg. {\bf 60} (1989),273-288.
\par
\medskip
\noindent \ref {L.2}: A.Lorenzini: {\it The minimal resolution conjecture.}
J. of  Alg. {\bf 156} (1993), 5-35.
\par
\medskip
\noindent \ref {PW}: P.Pragracz, J.Weyman: {\it Complexes associated with
trace and evaluation.} Adv. in Math. {\bf 57}, (1985), 163-207.
\par
\medskip
\noindent \ref {P-S}: C.Peskine, L.Szpiro: {\it Liaison de varietes
algebriques.} Inv. Math. {\bf 26}, (1974), 271-302.
\par
\medskip
\noindent \ref {R}: T.G.Room: {\it The Geometry of Determinantal Loci.}
Cambridge Univ. Press, {\bf 1938}, London.
\par
\medskip
\noindent \ref {S}: F.Steffen: {\it Generic determinantal
schemes and the smoothability of determinantal schemes of codimension 2.}
Manuscripta Math. {\bf 82}, (1994), 417-431.
\par
\medskip
\noindent \ref {S-T-V}: A.Simis, N.V.Trung, G.Valla: {\it The diagonal
subalgebra of a blow-up algebra.} Preprint.
\vskip 1cm
Addresses of the authors:
\bigskip
Alessandro Gimigliano: Dip. di Matematica, Univ. di Bologna,
\par
P.zza Porta S.Donato 5, I-40127; Bologna, Italy.
\par
e-mail:  GIMIGLIA@dm.unibo.it
\medskip
Anna Lorenzini: Dip. di Matematica, Univ. di Perugia,
\par
Via Vanvitelli 1, I-06123; Perugia, Italy.
\par
e-mail:  ANNALOR@ipguniv.unipg.it
\end